\documentclass[10pt]{amsart}
\usepackage{amssymb,latexsym,amsmath,amsfonts}
\usepackage[mathscr]{eucal}

\numberwithin{equation}{section}
\theoremstyle{definition}
\newtheorem{definition}{Definition}[section]

\theoremstyle{plain}
\newtheorem{theorem}[definition]{Theorem}
\newtheorem{lemma}[definition]{Lemma}

\newtheorem{example}[definition]{Example}

\newcommand{\CC}{\mathbb{C}^2}
\newcommand{\CD}{\mathbb{C}^m}
\newcommand{\cplx}{\mathbb{C}}
\newcommand{\RR}{\mathbb{R}^2}

\newcommand{\bdy}{\partial D}
\newcommand{\dom}{\mathscr{D}}
\newcommand{\DOM}{\widetilde{\mathscr{D_{\delta^*}}}}
\newcommand{\OM}{\Omega}
\newcommand{\ann}{{\rm Ann}}
\newcommand{\MO}{\widetilde{\Omega \ {}}}
\newcommand{\MOM}{\widetilde{\Omega_{1}}}
\newcommand{\poly}{\Delta_*}
\newcommand{\Dsc}{\boldsymbol{\sf P}}

\newcommand{\smoo}{\mathcal{C}}
\newcommand{\hol}{\mathcal{O}}

\newcommand{\gr}{\varGamma}
\newcommand{\diss}{\mathcal{D}}
\newcommand{\anal}{\mathscr{A}}
\newcommand{\banal}{\widetilde{\mathscr{A}_r}}
\newcommand{\biss}{\widetilde{\mathcal{D}_r}}
\newcommand{\pis}{\widetilde{\psi(}}
\newcommand{\crv}{\mathscr{C}}

\newcommand{\eps}{\varepsilon}

\newcommand{\zt}{\zeta}
\newcommand{\zbar}{\overline{z}}
\newcommand{\tht}{\theta}
\newcommand{\zahl}{\mathbb{Z}}

\newcommand{\exdt}{\widetilde{{\rm j} \ {}}}
\newcommand{\kup}{\mathfrak{K}}

\begin{document}

\title[Extension theorems of Hartogs-Chirka type]{The role of Fourier modes in\\
extension theorems of Hartogs-Chirka type}

\author{David E. Barrett}
\address{Department of Mathematics, University of Michigan, 525 East University Ann Arbor, MI
48109}
\email{barrett@umich.edu}
\thanks{The work of the first author is supported by NSF Grant DMS-0072237.}
\author{Gautam Bharali}
\email{bharali@umich.edu}

\begin{abstract} We generalize Chirka's theorem on the extension of functions
holomorphic in a neighbourhood of $\gr(F)\cup(\bdy\times D)$ -- where $D$ is
the open unit disc and $\gr(F)$ is the graph of a continuous $D$-valued function
$F$ -- to the bidisc. We extend holomorphic functions by applying the
Kontinuit{\"a}tssatz to certain continuous families of analytic annuli, which is
a procedure suited to configurations not covered by Chirka's theorem.
\end{abstract}                        
\maketitle

\section{Introduction and statement of results}\label{S:intro}

This article is motivated by the paper \cite{chirka:gHlnldbar?} by Chirka, in which
the following theorem is proved (in what follows, $D$ will represent the
open unit disc in $\cplx$ with centre at the origin, and given a function $F$ defined
in some region in $\cplx$, $\gr(F)$ will denote the graph of $F$ over its domain) :
\medskip

\noindent{{\bf Theorem (Chirka).} {\em Let $F\in\smoo(\overline{D};\cplx)$ and assume 
that $\sup_{\overline{D}}|F|<1$. Let $\OM$ be a connected neighbourhood of 
$\gr(F)\cup(\bdy\times D)$ contained in $\cplx\times D$. If $f\in\hol(\OM)$, then
$f$ extends holomorphically to the bidisc $D\times D$.}}
\smallskip

The requirement that $\sup_{\overline{D}}|F|<1$ is rather essential to the
extension theorem stated above (in contrast, refer to \cite{chirkaRosay:rpgHl98} for
a version by Chirka \& Rosay, in which the condition $\sup_{\overline{D}}|F|<1$ is 
relaxed, but in which only the functions holomorphic in the union of a neighbourhood of 
$\gr(F)$ with $\{z\in\cplx \ : \ |z|>1\}\times D$ -- i.e. holomorphic in a {\em large}
domain -- extend holomorphically). A pertinent
counterexample, when $\sup_{\overline{D}}|F|>1$, to the sort of holomorphic extension 
described in Chirka's theorem -- i.e. extension from {\em small} neighbourhoods of 
$\gr(F)\cup(\bdy\times D)$ -- is the case when $\gr(F)$ is a Wermer disc. We will 
discuss this example in \S \ref{S:ex} below. 
\medskip

The strategy of Chirka -- inspired by the methods in 
\cite{ivashkovichShevchishin:dncc&ems98} -- is to construct a continuously varying 
family of functions 
$\{F_t\}_{t\in[0,1]}\subset\{G\in\smoo(\cplx)  \ | \ \lim_{|z|\to\infty}G(z)=0\}$
such that $F_1=\widetilde{F \ {}}$ and $F_0\equiv 0$, and such that $\gr(F_t)$ is 
complex-analytic in a neighbourhood of any 
$(z,F_t(z))\notin\OM\cup(\{ \ |z|>1\}\times D)$. Here, $\widetilde{F \ {}}$ is any 
smooth extension of the $F$ provided by the theorem, that satisfies
$\widetilde{F|}_{|z|\geq 2}\equiv 0$. Next, one extends $f\in\hol(\OM)$ to
$\OM\cup(\{ \ |z|>1\}\times D)$ via Laurent decomposition. One can now show that the latter
can be analytically continued, owing to the Kontinuit{\"a}tssatz, via $\{F_t\}_{t\in[0,1]}$
to a neighbourhood of the classical Hartogs configuration 
$\gr(F_0)\cup(\bdy\times D)$. The condition $\sup_{\overline{D}}|F|<1$ is 
crucial in ensuring that the extension of $f\in\hol(\OM)$ by Laurent decomposition is 
single-valued. The strategy described fails in $\cplx^n, \ n >2$, and Chirka's theorem 
does not extend to higher dimensions as shown by Rosay's counterexample in \cite{rosay:crHp98}.
\medskip

The results in this paper are motivated by the two-fold aim of :
\begin{enumerate}
\item[a)] Showing that functions holomorphic in {\em small} neighbourhoods of a 
Hartogs-Chirka type configuration $\gr(F)\cup(\bdy\times D)$, with $\sup_{\overline{D}}|F|\gg 1$,
extend holomorphically to $D\times D$ (in a manner that will be made precise in Theorem 
\ref{T:mainC2}), given that $F$ satisfies suitable restrictions.
\item[b)] Extending Chirka's theorem to higher dimensions, and to a reasonably wide
class of Hartogs-Chirka type configurations $\gr(F)\cup(\bdy\times D^m), \ m\geq 2$ (in
particular, to configurations in which $F$ is not {\em merely} real-analytic or $\smoo^\infty$).
\end{enumerate}
Neither of the above seems to be achievable using Chirka's strategy. In this article, we discuss 
an alternative strategy for invoking the Kontinuit{\"a}tssatz, and use it to demonstrate
new Hartogs-Chirka type extension phenomena.
\medskip

The first of the above aims is met by the following theorem. But we first present the following
notation : if $\OM$ is a domain in $\cplx^n$, then 
$(\MO,\pi^{\OM})$ will denote the envelope of holomorphy of $\OM$.
\medskip

\begin{theorem}\label{T:mainC2} Let $F\in\smoo^(\overline{D};\cplx)$ and assume that
$\sup_{\bdy}|F|<1$. Let $A_j(r)$ represent the $j^{th}$ Fourier coefficient of 
$F(re^{i\centerdot})$, $r>0, \ j\in\zahl$. Assume that $F$ satisfies the condition
\begin{equation}\label{E:condtnC2}
\sum_{n\in\zahl}\frac{|A_n(r)|}{r^n} \ < \ 1 \ \forall r\in(0,1].
\end{equation}
Let $\OM_1$ be a neighbourhood of $\gr(F)\cup(\bdy\times D)$ and let $\OM_2$ be any
connected open set satisfying $\bdy\times D\subset\overline{\OM_2}
\subset\OM_1\cap(\{ \ |z|\geq 1\}\times D)$.
If $f\in\hol(\OM_1)$, then $f|_{\OM_2}$ has a holomorphic extension to $D\times D$.
\end{theorem}
\medskip

Note that since $\sup_{\overline{D}}|F|>1$, $\pi^{\OM_1}(\MOM)\varsupsetneq D\times D$
in general. For this reason, the usual arguments justifying that $f$ has a single-valued 
extension to the bidisc fail.
This is the reason behind the particular form of the conclusion of Theorem \ref{T:mainC2}.
Observe that while the condition \eqref{E:condtnC2} admits $F$ such that the negative
Fourier modes of $F(re^{i\centerdot})$ are large, it imposes a severe restriction
on the sizes of the positive Fourier modes of $F(re^{i\centerdot})$ as $r\to 0^+$. One
would like to investigate if such severe restrictions on the positive Fourier modes
are necessary. This is a valid concern because if we assume that the function $F$ has 
{\em only} positive Fourier modes, the condition \eqref{E:condtnC2} becomes 
unnecessary. The relevant theorem in this case is
\medskip

\begin{theorem}\label{T:subsC2} Let $F\in\smoo(\overline{D};\cplx)$ and assume that
$\sup_{\bdy}|F|<1$. Let $A_j(r)$ represent the $j^{th}$ Fourier coefficient of 
$F(re^{i\centerdot})$, $r>0, \ j\in\zahl$. Assume that $A_j\equiv 0 \ \forall j<0$.
Let $\OM_1$ be a neighbourhood of $\gr(F)\cup(\bdy\times D)$ and let $\OM_2$ be any
connected open set satisfying $\bdy\times D\subset\overline{\OM_2}
\subset\OM_1\cap(\{ \ |z|\geq 1\}\times D)$.
If $f\in\hol(\OM_1)$, then $f|_{\OM_2}$ has a holomorphic extension to $D\times D$.
\end{theorem}
\medskip

In a somewhat different direction, we may consider a 
continuous mapping $F:=(F_1,\dots,F_m) : \overline{D}\to D^m$, $m\geq 2$,
and consider the Hartogs-Chirka type configuration $\gr(F)\cup(\bdy\times D^m)$.
We know that, in general, Chirka's result is not true for such higher-dimensional 
configurations -- see \cite{rosay:crHp98}. In contrast, it has been shown by 
Bharali \cite{bharali:sgCet00} that Chirka's result does generalize to a certain class of
Hartogs-Chirka type configurations. However, the class of real-analytic maps 
$(F_1,\dots,F_m)$ studied in \cite{bharali:sgCet00} is rather restrictive. We show in 
this paper that that if we impose a 
condition analogous to condition \eqref{E:condtnC2} above, we can demonstrate analytic 
continuation for a considerably less restrictive set of configurations. We make this 
precise in the following 
\medskip
 
\begin{theorem}\label{T:mainCm} Let $F=(F_1,\dots,F_m)\in\smoo(\overline{D};\CD)$. Assume
that $F(e^{i\theta})\in D^m \ \forall\theta\in[0,2\pi)$ and let
$A_{jk}(r)$ represent the $k^{th}$ Fourier coefficient of $F_j(re^{i\centerdot})$, $r>0, \ 
k\in\zahl, \ j=1,\dots,m$. Assume that each $F_j$ satisfies the condition
\begin{equation}\label{E:condtnCm}
\sum_{n\in\zahl}\frac{|A_{jn}(r)|}{r^n} \ < \ 1 \ \forall r\in(0,1].
\end{equation}
Let $\OM_1$ be a neighbourhood of $\gr(F)\cup(\bdy\times D^m)$ and let $\OM_2$ be any
connected open set satisfying $\bdy\times D^m\subset\overline{\OM_2}
\subset\OM_1\cap(\{ \ |z|\geq 1\}\times D^m)$.
If $f\in\hol(\OM_1)$, then $f|_{\OM_2}$ has a holomorphic extension to $D\times D$.
\end{theorem}
\smallskip 

We note that if, in Theorem \ref{T:mainCm}, $F$ were to satisfy the restriction 
$F(\zt)\in D^m \ \forall\zt\in\overline{D}$, then all functions $f\in\hol(\OM)$ -- where
$\OM$ is a connected neighbourhood of $\gr(F)\cup(\bdy\times D^m)$ contained in
$\cplx\times D^m$ -- would extend to $D^m$, which is just Chirka's extension 
phenomenon in a restricted, higher-dimensional setting.
\medskip

The approach used in the first and the third theorem is to construct a continuous family
of analytic annuli which are attached to $\gr(H)$ -- where $H$ is an appropriately selected
perturbation of $F$ -- along their inner boundaries, and to
$\bdy\times D^m$ (with $m=1$ in Theorem \ref{T:mainC2} and $m\geq 2$ in Theorem 
\ref{T:mainCm}) along their outer boundaries. Once this family is constructed, analytic 
continuation is achieved by invoking the Kontinuit{\"a}tssatz. The proof of Theorem 
\ref{T:subsC2} uses a similar idea, but involves continuous families of analytic discs. 
These proofs may be found in \S \ref{S:proofs}. The technical construction of the 
aforementioned families of annuli/discs is carried out in the next section.
\medskip

In the final section of this paper, we discuss a few examples. Firstly, we show that one 
can construct Hartogs-Chirka type configurations $\gr(F)\cup(\bdy\times D)$ such that 
$\sup_{\overline{D}}|F|$ is as large as we want and such that functions holomorphic in 
small neighbourhoods of this configuration extend. Next, we discuss a configuration
involving Wermer's disc (see Example \ref{Ex:Wermer} for a definition) -- for which the 
extension phenomenon occuring in the previous example fails. And lastly, we show how 
Rosay's counterexample to a higher-dimensional analogue of Chirka's theorem fails to 
satisfy the hypotheses of Theorem \ref{T:mainCm}. 
\bigskip

\section{Preliminary lemmas}\label{S:vorspiel}

We need a few preliminary lemmas before we can prove our main theorems. In what follows,
$\ann(a;r,R)$ will denote the open annulus with centre at $a\in\cplx$ and having inner and
outer radii $r$ and $R$ respectively, while $D(a;R)$ will denote the open disc of
radius $R$ with centre at $a$. The symbol $\smoo^\infty(\overline{D};\cplx^m), \ m=1,2,\dots$,
will denote the class of infinitely differentiable functions on the unit disc, all of whose
derivatives extend to continuous functions on $\overline{D}$.
\medskip

The reader will notice that in the following lemma the hypothesis 
$G\in\smoo^\infty(\overline{D};\cplx)$ is much stronger than is required for the conclusion
of Lemma \ref{L:1stlem}. The only place where we use this hypothesis is in showing the 
existence of a certain limit towards the end of the proof. However, stating the strongest
versions of Lemmas \ref{L:1stlem} and \ref{L:2ndlem} -- which are of relatively minimal
utility in themselves -- merely results in statements that are overly technical. For this
reason, the $G$ occuring in Lemmas \ref{L:1stlem}-\ref{L:dhull} shall be assumed to be 
$\smoo^\infty$.
\smallskip

\begin{lemma}\label{L:1stlem} Let $G(re^{i\tht})=\sum_{n=-N}^Nb_n(r)e^{in\tht}$ and assume
that $G\in\smoo^\infty(\overline{D};\cplx)$. Assume further that
\begin{equation}\label{E:1/contact}
\sum_{n=-N}^N\frac{|b_n(r)|}{r^n} \ < \ 1 \ \forall r\in(0,1].
\end{equation}
Then the holomorphic function
\[
\anal_r(\zt) \ = \sum_{n=-N}^N b_n(r)\left(\frac{\zt}{r}\right)^n, \quad\zt\in\ann(0;r,1),
\]
which belongs to $\hol[\ann(0;r,1)]\cap\smoo[\overline{\ann}(0;r,1)]$, satisfies
$|\anal_r(e^{i\theta})|<1$. Fix $\nu\in\mathbb{N}$ and let $K\Subset\ann(0;1/\nu,1)$ be a 
compact subset. The function $(0,1/\nu]\times K\ni (r,\zt)\mapsto\anal_r(\zt)$ 
extends to a continuous function on $[0,1/\nu]\times K$.
\end{lemma}
\begin{proof}
To prove the first part of this lemma, note that
\begin{equation}\label{E:conv}
|\anal_r(e^{i\theta})| \ \leq \ \sum_{n=-N}^N|b_n(r)|\left|\frac{e^{i\theta}}{r}\right|^n \ = \
	\sum_{n=-N}^N\frac{|b_n(r)|}{r^n} \ < \ 1.
\end{equation}
We fix $\nu\in\mathbb{N}$ and then fix a compact set $K\Subset\ann(0;1/\nu,1)$.
It is obvious that 
$(0,\nu]\times K\ni (r,\zt)\mapsto b_{-n}(r)(r/\zt)^n$ extends to a continuous function
on $[0,1/\nu]\times K$, which simply vanishes when $r=0$, for each $n=1,2,\dots,N$. 
Now consider the function 
$(r,\zt)\mapsto b_{n}(r)(\zt/r)^n, \ n=1,2,\dots,N$. Note that \eqref{E:1/contact} 
$\Longrightarrow \ |b_n(r)|<r^n \ \forall n=1,2,\dots,N$. This implies, since $G$ is assumed 
to be smooth, that each of the latter functions extends continuously to a function
$\varphi_n\in\smoo([0,1/\nu]\times K)$, which is defined as
\[
\varphi_n(r,\zt) \ := \ \begin{cases}
			b_{n}(r)(\zt/r)^n, &\text{if $(r,\zt)\in(0,1/\nu]\times K$},\\
			{} & {} \\
			\dfrac{1}{n!} \left.\dfrac{d^nb_n}{dr^n}\right|_{r=0}\zt^n,
			&\text{if $(r,\zt)\in\{0\}\times K$.}
\end{cases}
\]
Since $\anal_r$ is a finite sum of the functions $b_{n}(r)(r/\zt)^n$, the last two 
observations establish the second part of this lemma.
\end{proof}

\begin{lemma}\label{L:hull} Let $G$ be as in Lemma \ref{L:1stlem}, but assume additionally that
$\sup_{\bdy}|G|<1$. Let $\OM_1$ be a neighbourhood of $\gr(G)\cup(\bdy\times D)$ such as that
described in Theorem \ref{T:mainC2}. Then
\begin{enumerate}
\item[a)] $\{\anal_r\}_{r\in(0,1)}$ is a continuous family in the sense that for a fixed 
$\zt_0\in D\setminus\{0\}$, $r\mapsto\anal_r(\zt_0)$ is continuous in the interval
$(0,|\zt_0| \ )$.
\item[b)] $\lim_{r\to 0^+}\anal_r(\zt_0)$ exists for each $\zt_0\in D\setminus\{0\}$,
and there exists a $\psi\in\hol(D)$ such that 
$\psi(\zt)=\lim_{r\to 0^+}\anal_r(\zt)$ on $D\setminus\{0\}$.
\item[c)] Define 
\[
\kup \ := \ \gr(\psi)\cup\left[ \ \cup_{0<r<1}\{(\zt,\anal_r(\zt))\in\CC \ | \ r<|\zt|<1\} \ 
		\right]\setminus\OM_1.
\]
$\kup$ is compact.
\end{enumerate}
\end{lemma}
\begin{proof}
Part (a) and the first half of part (b) are obvious conclusions of Lemma \ref{L:1stlem}.
Thus, we may define
\[
\psi(\zt) \ := \ \lim_{|\zt|>r\to 0^+}\anal_r(\zt) \quad\forall\zt\in D\setminus\{0\}.
\]
Fix a $\nu\in\mathbb{N}$. Lemma \ref{L:1stlem}(b) tells us that
\[
(\anal_r|_{\ann(0;1/\nu,1)})(\zt) \ \longrightarrow \psi(\zt) \ 
\text{\em uniformly on each compact $K\Subset\ann(0;1/\nu,1)$ as $r\searrow 0$ .}
\]
We conclude from this statement that 
\begin{equation}\label{E:holom}
\psi|_{\ann(0;1/\nu,1)} \ \in \ \hol[\ann(0;1/\nu,1)] \quad \forall\nu=2,3,4,\dots
\end{equation}

Before proceeding any further, we comment that the functions $\anal_r$ are so constructed that 
$\gr(\anal_r)$, $0<r<1$, are analytic annuli that are attached to $\gr(G)$ along their inner 
boundaries and -- in view of the inequality \eqref{E:conv}  -- to $\bdy\times D$ 
along their outer boundaries. Therefore, 
\begin{multline}\label{E:upbd}
|\anal_r(\zt)| \ \leq \ 
\max\left\{\sup_{|\xi|=r}|\anal_r(\xi)|, \ 1\right\} \ \leq \ 
\max\left\{\sup_{\overline{D}}|G|, \ 1\right\} \\
	\forall \zt\in\ann(0;r,1) \ \text{\em and for each $r\in(0,1)$.}
\end{multline}
By \eqref{E:holom}, $\psi$ is already holomorphic on $D\setminus\{0\}$.
The bounds above imply, since $\psi(\zt)$ is the limit of the $\anal_r(\zt)$'s, provided
$\zt\neq 0$, that $|\psi(\zt)|\leq\sup_{\xi\neq 0}|G(\xi)|$ in a punctured neighbourhood of 
the origin. Thus, $\psi$ extends to a holomorphic function on $D$. This establishes (b).
\medskip

Notice that by the estimates \eqref{E:upbd} and part (b) of this lemma, $\kup$ is a bounded
set. Therefore, it suffices to show that $\kup$ is closed. Now consider a point 
$(z,w)\notin\OM_1$ with the property that there exist sequences
$\{r(\nu)\}_{\nu\in\mathbb{N}}\subset(0,1)$ and $\{\zt_\nu\}_{\nu\in\mathbb{N}}\subset D$
such that $r(\nu)\to 0$ as $\nu\to\infty$ and 
$(\zt_\nu,\anal_{r(\nu)}(\zt_\nu))\longrightarrow (z,w)$ as $\nu\to\infty$. 
It is easy to see that to prove (c),
it suffices to show that all such points $(z,w)\in\kup$. Notice that, by construction,
there is a $\delta(\OM_1)>0$ depending only on $\OM_1$ such that 
$(\zt,\anal_r(\zt))\in\OM_1 \ \forall r,|\zt|<\delta(\OM_1)$. Thus,
as $(z,w)\notin\OM_1$, $z\neq 0$. Now, given that the $\anal_r$'s converge uniformly
on compact subsets lying away from $0$, there exists $\kappa_1\in\mathbb{N}$ such that
\[
|\anal_{r(\nu)}(\zt)-\psi(\zt)| \ < \ \eps/2 \quad\forall \nu\geq \kappa_1, \ 
						\forall\zt\in\overline{D}(z;|z|/2).
\]
Let $\kappa_2\in\mathbb{N}$ be such that 
\[
\zt_\nu\in D(z;|z|/2) \qquad\text{\em and}\qquad
|w-\anal_{r(\nu)}(\zt_\nu)| \ < \ \eps/2 \quad \forall\nu\geq \kappa_2.
\]
The above inequalities imply that
\[
|w-\psi(\zt_\nu)| \ \leq \ |w-\anal_{r(\nu)}(\zt_\nu)|+
		|\anal_{r(\nu)}(\zt_\nu)-\psi(\zt_\nu)| \ < \eps \quad
		\forall \nu\geq\max(\kappa_1,\kappa_2).
\]
This tells us that $w=\lim_{\nu\to\infty}\psi(\zt_\nu)$, whence $(z,w)\in\gr(\psi)\setminus\OM_1$.
This establishes (c), and concludes our proof.
\end{proof}
\medskip

\begin{lemma}\label{L:2ndlem} Let $G(re^{i\tht})=\sum_{n=0}^Nb_n(r)e^{in\tht}$ -- i.e. we assume
that $G(re^{i\centerdot})$ has no negative Fourier modes. Assume further
that $G\in\smoo^\infty(\overline{D};\cplx)$. 
Then the holomorphic function
\[
\diss_r(\zt) \ = \sum_{n=0}^N b_n(r)\left(\frac{\zt}{r}\right)^n, \quad\zt\in D(0;r),
\]
which belongs to $\hol[D(0;r)]\cap\smoo[\overline{D}(0;r)]$, satisfies 
$\diss_r(re^{i\theta})=G(re^{i\theta}) \ \forall\theta\in[0,2\pi)$.
Fix $\nu\in\mathbb{N}$ and let $K\Subset D(0;1-1/\nu)$ be a 
compact subset. The function $(r,\zt)\mapsto\diss_r(\zt)$ 
is a continuous function on $[1-1/\nu,1]\times K$.
\end{lemma}
\smallskip
The above lemma is a triviality; we merely state it as an element that will be needed in the 
proof of our next result.
\medskip

\begin{lemma}\label{L:dhull} Let $G$ be as in Lemma \ref{L:2ndlem}, but assume additionally that
$\sup_{\bdy}|G|<1$. Let $\OM_1$ be a neighbourhood of $\gr(G)\cup(\bdy\times D)$ such as that
described in Theorem \ref{T:subsC2}. Then
\begin{enumerate}
\item[a)] $\{\diss_r\}_{r\in(0,1)}$ is a continuous family in the sense that for a fixed 
$\zt_0\in D$, $r\mapsto\diss_r(\zt_0)$ is continuous in the interval
$( \ |\zt_0|, 1)$.
\item[b)] $\lim_{r\to 1^-}\diss_r(\zt)$ exists for each $\zt\in D$, and this
limit defines a holomorphic function $\psi\in\hol(D)$.
\item[c)] Define 
\[
\kup \ := \ \gr(\psi)\cup\left[ \ \cup_{0<r<1}\{(\zt,\diss_r(\zt))\in\CC \ | \ |\zt|<r\} \ 
		\right]\setminus\OM_1.
\]
$\kup$ is compact.
\end{enumerate}
\end{lemma}
\begin{proof}
Part (a) and the first half of part (b) are direct consequences of Lemma \ref{L:2ndlem}.
The inference that 
\[
\psi(\zt) \ := \ \lim_{|\zt|<r\to 1^-}\diss_r(\zt) \quad\forall\zt\in D
\]
is holomorphic follows from Lemma \ref{L:2ndlem}. The uniform-convergence argument is
exactly analogous to the argument used in proving Lemma \ref{L:hull}. We therefore
omit the details. We remark that 
\[
\psi(\zt) \ = \ \sum_{n=0}^N b_n(1)\zt^{n}.
\]
\smallskip

The functions $\diss_r$ are so constructed that $\gr(\diss_r)$,
$0<r<1$, are analytic discs that are attached to $\gr(G)$ along their boundaries.
Therefore, 
\[
|\diss_r(\zt)| \ \leq \ 
\sup_{|\xi|=r}|\diss_r(\xi)| \ \leq \ \sup_{\overline{D}}|G|
\quad\forall\zt\in D(0;r) \ \text{\em and for each $r\geq 1-1/\nu$.}
\]
Thus, $\kup$ is a bounded set, and we argue that $\kup$ is closed exactly as we
did in Lemma \ref{L:hull}(c).
\end{proof}
\medskip

The following lemma is key to the proofs of Theorems \ref{T:mainC2}-\ref{T:mainCm}. Before
proving it, we explicitly state the following simple
\smallskip

\noindent{{\bf Fact :} {\em Due to the continuity of the functions $F$ and $G$ occuring in
the statements of the various theorems and lemmas above, the associated Fourier coefficients
$A_n(r)$ and $b_n(r)$ satisfy $A_n(0)=b_n(0)=0 \ \forall n\neq 0$.}
\smallskip

\noindent{This fact is used implicitly at several places in the next lemma.}
\smallskip

\begin{lemma}\label{L:weird} Let $F\in\smoo(\overline{D};\cplx)$ and let $A_j(r)$ 
represent the $j^{th}$ Fourier coefficient of $F(re^{i\centerdot})$, $r>0, \ j\in\zahl$. 
Assume that :
\begin{enumerate}
\item[1)] $\sup_{\bdy}|F|<1$, and
\item[2)] $F$ satisfies the condition
\begin{equation}\label{E:moments}
\sum_{n\in\zahl}\frac{|A_n(r)|}{r^n} \ < \ 1 \quad \forall r\in(0,1].
\end{equation}
\end{enumerate}
Given $\eps>0$ there exists a function $G\in\smoo^\infty(\overline{D};\cplx)$ of the form
\[
G(re^{i\theta}) \ = \ \sum_{n=-N}^N B_n(r)e^{in\theta},
\] 
where $N$ is some large positive integer and $B_n\in\smoo^\infty([0,1];\cplx)$, such that
\begin{itemize}
\item $|F(\zt)-G(\zt)|<\eps \ \forall\zt\in\overline{D}$,
\item $G$ has the property (1) and satisfies the analogue of (2) above 
(with $B_n(r)$ replacing $A_n(r)$ in \eqref{E:moments} above).
\end{itemize}
Furthermore, if property (2) is replaced by
\begin{itemize}
\item[$2^*$)] $F$ has no negative Fourier modes, 
\end{itemize}
then $G$ can be constructed so that it has property (1) and $B_{-j}\equiv 0$ for $j=1,2,\dots,N$.
\end{lemma}
\begin{proof}
Define 
\begin{align}
S_m(\theta,r) \ &:= \ \sum_{j=-m}^mA_j(r)e^{ij\theta}\notag \\
\sigma_n(\theta,r) \ &:= \ \frac{S_0(r,\theta)+\dots+S_n(\theta,r)}{n+1}\notag
\end{align}
Let us first assume that $F$ has properties (1) and (2).
Let $\eta>0$ be so small that
\begin{align}\label{E:border}
\eta \ &< \ 1-\sup_{\bdy}|F|,\\
\eta+\sum_{n\in\zahl}|A_n(r)|/r^n \ &< \ 1 \quad\ \forall r\in(0;1],\notag
\end{align}
and define $\delta:=\min(\eps,\eta)$. There exists a natural number $N>0$ such that
\begin{equation}\label{E:Fejer}
|F(re^{i\theta})-\sigma_N(\theta,r)| \ < \delta/2 \quad \forall(\theta,r)\in[0,2\pi)\times[0,1].
\end{equation}
This above is a consequence of Fej{\'e}r's theorem. For a {\em fixed} $r\in[0,1]$, 
\eqref{E:Fejer} is precisely the statement of Fej{\'e}r's theorem applied to the
periodic function $F(re^{i\centerdot})$. However, on examining the proof of Fej{\'e}r's
theorem, one sees that owing to the equicontinuity of the family 
$\{F(re^{i\centerdot})\}_{r\in[0,1]}\subset\smoo(\mathbb{T})$, the
choice of $N$ in \eqref{E:Fejer} is uniform in $r\in[0,1]$.
\medskip

One sees immediately that if one writes 
\[
\sigma_N(\theta,r) \ = \ \sum_{j=-N}^{N}a_j(r)e^{ij\theta},
\]
then $|a_j(r)|\leq|A_j(r)| \ \forall r\in[0,1]$. 
For each $j=1,2,\dots,N$, we pick a function $B_{-j}(r)$ which satisfies 
the following conditions :
\begin{enumerate}
\item[$i$)] $B_{-j}\in\smoo^\infty([0,1];\cplx)$,
\item[$ii$)] $B_{-j}$ vanishes to infinite order at $r=0$, and
\item[$iii$)] $|a_{-j}(r)-B_{-j}(r)|\leq\delta/2(2N+1) \ \forall r\in[0,1]$,
\end{enumerate}
provided $a_{-j}\not\equiv 0$, $j=1,2,\dots,N$.  {\em If $a_{-j}\equiv 0$, we just
choose $B_{-j}\equiv 0$.}
Note that by our condition on the Fourier coefficients $\{A_n(r)\}_{n\in\zahl} \ $, 
$|a_j(r)| \leq|A_j(r)| \leq r^j \ \forall j=1,2,\dots,N, \ r\in[0,1)$.
Let $0<R_0<1$ be a small number such that
\[
R_0 \ \leq \ \left\{\frac{\delta}{4(2N+1)}\right\}^{1/j}\quad \forall j=1,2,\dots,N.
\] 
For each $j=1,2,\dots,N$ we define a function $B_j(r)$ as follows :
\[
B_j(r) \ := \ \begin{cases}
		\alpha_j(r)r^j, &\text{if $r\leq R_0$},\\
		\beta_j(r), &\text{if $r\geq R_0$},
		\end{cases}
\]
such that
\begin{enumerate}
\item[$i^*$)] $B_j\in\smoo^\infty([0,1];\cplx)$,
\item[$ii^*$)] $\alpha_j$ vanishes to infinite order at $r=0$,
\item[$iii^*$)] $\alpha_j$ satisfies
\[
|\alpha_j(r)| \ \leq \ \sup_{s\leq 1}\frac{|a_j(s)|}{s^j} \quad\forall s\in[0,R_0],
\]
\item[$iv^*$)] $\beta_j$ satisfies
\[
|\beta_j(r)-a_j(r)| \ \leq \ \frac{R_0^j\delta}{2(2N+1)} \quad\forall r\in[R_0,1].
\]
\end{enumerate}
Finally, define $B_0(r)$ to be any $\smoo^\infty$ function such that 
$|B_0(r)-a_0(r)|<\delta/2(2N+1) \ \forall r\in[0,1]$ and such that 
$B_0-B_0(0)$ vanishes to high order at $r=0$. Now write
\[
G(re^{i\theta}) \ = \ \sum_{j=-N}^N B_j(r)e^{ij\theta}.
\]
\smallskip

We now make some estimates. We first consider the $B_{-j}(r)$'s, $j=1,2,\dots,N$.
Note that the following statements continue to be true trivially if $B_{-j}\equiv 0$
for any $j=1,2,\dots,N$. 
\begin{align}\label{E:negbeer}
\sum_{j=1}^{N}|B_{-j}(r)-a_{-j}(r)| \ &\leq \ \sum_{j=1}^N\frac{\delta}{2(2N+1)} \
		\leq \ \frac{\delta N}{2(2N+1)}, \\
\sum_{j=1}^{N}|B_{-j}(r)|r^j \ &\leq \ \sum_{j=1}^N r^j\left\{|a_{-j}(r)|+
						\frac{\delta}{2(2N+1)}\right\}\label{E:negbee}\\
		\ &\leq \ \sum_{j=1}^N|a_{-j}(r)|r^j + \frac{\delta N}{2(2N+1)} 
				\quad\forall r\in[0,1].\notag
\end{align}
Next, we consider the $B_{j}(r)$'s, $j=1,2,\dots,N$. First, we let $0\leq r\leq R_0$. We
use item ($iii^*$) in the definition of $B_j(r)$ above to get : 
\begin{align}\label{E:posbeer1}
\sum_{j=1}^{N}|B_{j}(r)-a_j(r)| \ &\leq \ \sum_{j=1}^N r^j\left|\alpha_j(r)-
					\frac{|a_j(r)|}{r^j}\right| \\
	&\leq \ \sum_{j=1}^N 2R_0^j \ \sup_{s\leq 1}\frac{|a_j(s)|}{s^j}\notag\\
	&\leq \ \sum_{j=1}^N 2\left\{\frac{\delta}{4(2N+1)}\right\} \ \sup_{s\leq 1}
		\frac{|A_j(s)|}{s^j} \notag\\
	& \leq \ \frac{\delta N}{2(2N+1)} \ , \notag\\
\sum_{j=1}^{N}\frac{|B_{j}(r)|}{r^j} \ &\leq \ \sum_{j=1}^N \ 
		\sup_{s\leq 1}\frac{|a_{j}(s)|}{s^j} \quad\forall r\in[0,R_0].\label{E:posbee1}
\end{align}
And when we consider $R_0\leq r\leq 1$, we use item ($iv^*$) in the definition of $B_j(r)$
to get :
\begin{align}\label{E:posbeer2}
\sum_{j=1}^{N}|B_{j}(r)-a_j(r)| \ &\leq \ \sum_{j=1}^N\frac{R_0^j\delta}{2(2N+1)} \ 
		\leq \ \frac{\delta N}{2(2N+1)}, \\
\sum_{j=1}^{N}\frac{|B_{j}(r)|}{r^j} \ &\leq \ \sum_{j=1}^N 
\left\{\frac{|a_{j}(r)|}{r^j}+\frac{R_0^j\delta}{2(2N+1)R_0^j}\right\} \label{E:posbee2}\\
		\ &\leq \ \sum_{j=1}^N\frac{|a_{j}(r)|}{r^j} + \frac{\delta N}{2(2N+1)} 
				\quad\forall r\in[R_0,1].\notag
\end{align}

Observe that from \eqref{E:negbee} and \eqref{E:posbee1}, we have
\begin{align}
\sum_{j=-N}^{N}\frac{|B_{j}(r)|}{r^j} 
\ &\leq \ \sum_{j=1}^N|a_{-j}(r)|r^j + \frac{\delta N}{2(2N+1)} + |a_0(r)| + 
		\frac{\delta}{2(2N+1)} + \sum_{j=1}^N \ \sup_{s\leq 1}\frac{|a_{j}(s)|}{s^j} \notag\\
	&\leq \ \sum_{j=-N}^N \ \sup_{s\leq 1}\frac{|a_{j}(s)|}{s^j} + 
		\frac{\delta(N+1)}{2(2N+1)} \notag\\
	&\leq \ \sum_{j=-N}^N \ \sup_{s\leq 1}\frac{|A_{j}(s)|}{s^j} + 
		\frac{\eta}{2} \ < 1 \quad\forall r\in[0,R_0]. \label{E:posbird1}
\end{align}
The last inequality follows from the definition of $\eta$ and the fact that
$|a_j(r)|\leq|A_j(r)| \ \forall r\in[0,1]$. Next, applying \eqref{E:negbee} and
\eqref{E:posbee2} we get
\begin{align}
\sum_{j=-N}^{N}\frac{|B_{j}(r)|}{r^j} 
\ &\leq \ \sum_{j\neq 0}\frac{|a_{j}(r)|}{r^j} + \frac{2\delta N}{2(2N+1)} + |a_0(r)| + 
		\frac{\delta}{2(2N+1)} \notag\\
	&\leq \ \sum_{j=-N}^N \ \sup_{s\leq 1}\frac{|a_{j}(s)|}{s^j} + 
		\frac{\delta}{2} \notag\\
	&\leq \ \sum_{j=-N}^N \ \sup_{s\leq 1}\frac{|A_{j}(s)|}{s^j} + 
		\frac{\eta}{2} \ < 1 \quad\forall r\in[R_0,1]. \label{E:posbird2}
\end{align}

From the inequalities \eqref{E:border}, \eqref{E:posbird1} and \eqref{E:posbird2}, we get
\[
\sum_{j=-N}^{N}\frac{|B_{j}(r)|}{r^j} \ < \ 1 \quad\forall r\in[0,1],
\]
which is to say that $G$ satisfies the analogue of (2), with $B_n(r)$ replacing $A_n(r)$ in the
expression \eqref{E:moments}.
\medskip

We now exploit the estimates \eqref{E:negbeer}, \eqref{E:posbeer1} and \eqref{E:posbeer2},
to get
\begin{align}
|F(re^{i\theta})-G(re^{i\theta})| \ &\leq \ |F(re^{i\theta})-\sigma_N(re^{i\theta})| +
		\sum_{j=-N}^N|a_j(r)-B_j(r)|\notag\\
&< \ \frac{\delta}{2}+ 2\cdot\frac{\delta N}{2(2N+1)} + \frac{\delta}{2(2N+1)} \ = \delta.
\end{align}
Given the way in which $\delta$ is defined, we see that $G$ has the property (1), and 
$|G(\zt)-F(\zt)|<\eps \ \forall\zt\in\overline{D}$. $G$ is, of course, smooth by 
construction.
\medskip

Note further that in the above construction, if $F$ has no negative Fourier modes, neither 
does $G$. So, if $F$ had the property ($2^*$) instead of property (2), in addition to choosing
$B_{-j}\equiv 0$ we would use the
same rule by which we selected $B_0(r)$ in the above argument to define 
$B_j(r), \ j=1,2,\dots,N$. It is easy to verify that this modified construction would yield the 
second part of this lemma.
\end{proof}

\section{Proofs of the theorems}\label{S:proofs}

\begin{proof}[{\bf 3.1. The proof of Theorem \ref{T:mainC2}}] Let $\eps>0$ be so 
small that $F(D)+\eta\subset\Omega_1 \ \forall\eta\in\cplx$ such that $|\eta|<2\eps$.
By Lemma \ref{L:weird}, there exists a function $H\in\smoo^\infty(\overline{D};\cplx)$ 
with $H(re^{i\theta})=\sum_{n=-N}^N B_n(r)e^{in\theta}$ such that 
\begin{align}
|H(\zt)-F(\zt)| \ < \ \eps \quad\forall\zt\in\overline{D},\qquad
		\sup_{\bdy}|H| \ < \ 1,\notag\\
\sum_{n=-N}^N\frac{|B_n(r)|}{r^n} \ < \ 1 \quad\forall r\in(0,1].\notag
\end{align}
Let $\delta>0$ be so small that
\begin{itemize}
\item $\sup_{\bdy}|H|+\delta<1$;
\item $\delta+\sum_{n=-N}^N|B_n(r)|/r^n < 1 \ \forall r\in(0;1]$; and
\item $H(D)+\eta\subset\OM_1 \ \forall\eta\in\cplx$ such that $|\eta|<\delta$.
\end{itemize}
Define, for each $\eta$ such that $|\eta|<\delta$
\begin{align}
H^{(\eta)}(\zt) \ &:= \ H(\zt)+\eta, \notag \\
\anal_r^{(\eta)}(\zt) \ &:= \ \sum_{n\neq 0}B_n(r)\left(\frac{\zt}{r}\right)^n
					+(B_0(r)+\eta), \quad\zt\in\ann(0;r,1),\notag
\end{align}
We apply Lemma \ref{L:hull} to $\{\anal_r^{(\eta)}\}_{r\in(0;1)}$ for each $\eta$, by
taking
\[
b_n(r) \ = \ B_n(r) \quad\forall n\neq 0, \qquad b_0(r) \ = \ B_0(r)+\eta
\]
in that lemma. We conclude that there is a function $\psi\in\hol(D)$ such that
\begin{enumerate}
\item[1)] For any fixed $\zt_0\in D\setminus\{0\}$, $r\mapsto\anal_r^{(\eta)}(\zt_0)$ is continuous
in $(0,|\zt_0| \ )$ for $|\eta|<\delta$.
\item[2)] For $\zt\in D\setminus\{0\}$, $\lim_{r\to 0^+}\anal_r^{(\eta)}(\zt)=
\psi(\zt)+\eta$ for $|\eta|<\delta$.
\item[3)] For each $\eta \ : \ |\eta|<\delta$, $\kup^{(\eta)}$ is compact, where we define
\[
\kup^{(\eta)} \ := \ \gr(\psi+\eta)\cup
			\left[ \ \cup_{0<r<1}\{(\zt,\anal_r^{(\eta)}(\zt))\in\CC \ | 
			\ r<|\zt|<1\} \ \right]\setminus\OM_1.
\]
\end{enumerate}
In other words, for each fixed $\eta$, the family $\{\gr(\anal_r^{(\eta)})\}_{r\in(0,1)}$ is a 
continuous family of analytic annuli attached to $\gr(H^{(\eta)})\cup(\bdy\times D)$, which
accumulate onto $\gr(\psi+\eta)$ as $r\to 0^+$. The analytic annuli $\gr(\anal_r^{(\eta)})$
with $r\approx 1$ are contained in $\OM_1$. We can therefore apply the Kontinuit{\"a}tssatz to
conclude that 
\begin{equation}\label{E:nbd}
U(\Delta) \ := \ \{(z,w)\in D\times\cplx \ : \ (w-\psi(z))\in\Delta\} \ = \ 
\bigcup_{\eta\in\Delta}\gr(\psi+\eta)  \ \subset \ \pi^{\OM_1}(\MOM),
\end{equation}
where $\Delta$ is any disc contained in $D(0;\delta)$
\medskip

There is a canonical holomorphic imbedding of $\OM_1$
into $\MOM$. We denote this imbedding by ${\rm j}:\OM_1\hookrightarrow\MOM$. 
Corresponding to each $f\in\hol(\OM_1)$ there is a holomorphic function on $\MOM$,
which we shall denote by $\mathcal{E}(f)$, such that $\mathcal{E}(f)\circ{\rm j}=f$. 
It is now a standard argument -- see, for instance \cite{chirkaStout:K95}
or \cite{chirka:gHlnldbar?} -- to show
that there exist holomorphic mappings
\[
\banal(\centerdot \ ;\eta):\ann(0;r,1)\to\MOM \ \forall r\in(0,1), \ \text{and}\quad 
\pis\centerdot \ ;\eta):D\to\MOM
\]
such that
\begin{enumerate}
\item[a)] For each $\eta$ with $|\eta|<\delta$ :
\begin{align}
\pi^{\OM_1}\circ \banal(\zt;\eta)&=(\zt,\anal_r^{(\eta)}(\zt)) \ 
\forall\zt\in\ann(0;r,1) \ \text{when} \ r\in(0,1), \ \text{and}\notag\\ 
\pi^{\OM_1}\circ \pis\zt;\eta)&=(\zt,\psi(\zt)+\eta) \ \forall\zt\in D.\notag
\end{align}
\item[b)] ${\rm j}(\zt,\anal_r^{(\eta)}(\zt))=\banal(\zt;\eta)$ wherever the left-hand
side is defined, and $\forall r\in(0,1)$.
\item[c)] ${\rm j}(\zt,\psi(\zt)+\eta)=\pis\zt;\eta)$ wherever the left-hand
side is defined.
\end{enumerate}
Notice that -- in view of item (a) -- for each fixed $\zt\in D$, $\eta\mapsto 
\pis\zt;\eta)$
is holomorphic.
\medskip

Let $V$ be the connected component of $(\pi^{\OM_1})^{-1}(U( \ |\eta|<\delta))$ containing 
$\crv_0:={\rm image}(\pis\centerdot \ ;0))$. For each point $q\in\crv_0$ there 
is a neighbourhood
$W(q)\Subset V$ of $q$ such that $\pi^{\OM_1}|_{W(q)} : W(q)\to \CC$ is a biholomorphism.
Let $\poly$ be a disc centered at the origin that is so small that
\[
{\rm image}(\pis\centerdot \ ;\eta)) \ \subset \ \bigcup_{q\in\crv_0}W(q)\quad 
\forall\eta\in\poly.
\]
We define $\OM^*:=U(\poly)\cup\omega_2$, $\omega_2$ being a connected open set
satisfying 
\begin{itemize}
\item $\overline{\OM_2} \ \subset \ \omega_2 \ \subset \ \OM_1$; and
\item $\omega_2\cap U(\poly)$ is connected,
\end{itemize}
where $\OM_2$ is as described in Theorem \ref{T:mainC2},
and $U(\poly)$ is as defined by \eqref{E:nbd}. Our goal is to map $\OM^*$ into $\MOM$ in
such a way that this mapping extends ${\rm j}$. This mapping will allow us -- given any 
$f\in\hol(\OM_1)$ -- to extend $f|_{\omega_2}$ to $\OM^*$. But $\OM^*$ is a neighbourhood of 
a classical Hartogs configuration, whence $f|_{\OM_2}$ would extend to the bidisc.
To this end, we define
\[
\exdt(z,w) \ := \ \begin{cases}
			\pis z;w-\psi(z)), &\text{if $(z,w)\in U(\poly)$},\\
			{} &{}\\
			{\rm j}(z,w), &\text{if $(z,w)\in\omega_2$}.
			\end{cases}
\]
Note that if $(z,w)\in U(\poly)\cap\omega_2$, then, in view of item (c) above
\begin{equation}\label{E:coincide}
\pis z;w-\psi(z)) \ = \ {\rm j}(z,\psi(z)+\{w-\psi(z)\}) \ = \ {\rm j}(z,w)\quad 
\forall(z,w)\in U(\poly)\cap\omega_2.
\end{equation}
Thus, $\exdt$ is well-defined, and extends ${\rm j}$. From our foregoing remarks, 
$\exdt$ is holomorphic. Given any $f\in\hol(\OM_1)$, we define 
$\widetilde{f}\in\hol(\OM_*)$ by $\widetilde{f}:=\mathcal{E}(f)\circ\exdt$. In view
of \eqref{E:coincide}, $\widetilde{f}|_{\omega_2}\equiv f|_{\omega_2}$. Notice that
$\OM_*$ is a neighbourhood of $\gr(\psi)\cup(\bdy\times D)$, which is the classical
Hartogs configuration. Thus $\widetilde{f}$ has a holomorphic extension to $D\times D$, 
whence $f|_{\OM_2}$ has a holomorphic extension to $D\times D$. This concludes our
proof.
\end{proof}
\medskip

\begin{proof}[{\bf 3.2. The proof of Theorem \ref{T:subsC2}}] Since the proof of this 
theorem is similar to that of Theorem \ref{T:mainC2}, we shall be brief. 
Let $\eps>0$ be so 
small that $F(D)+\eta\subset\Omega_1 \ \forall\eta\in\cplx$ such that $|\eta|<2\eps$.
By Lemma \ref{L:weird}, there exists a function $H\in\smoo^\infty(\overline{D};\cplx)$ 
with $H(re^{i\theta})=\sum_{n=0}^N B_n(r)e^{in\theta}$ such that 
\[
|H(\zt)-F(\zt)| \ < \ \eps \quad\forall\zt\in\overline{D},\qquad
		\sup_{\bdy}|H| \ < \ 1.
\]
Note that the $H$ that has the two properties above can be so chosen that it has no
negative Fourier modes.  
Let $\delta>0$ be so small that
\begin{itemize}
\item $\sup_{\bdy}|H|+\delta<1$; and
\item $H(D)+\eta\subset\OM_1 \ \forall\eta\in\cplx$ such that $|\eta|<\delta$.
\end{itemize}
For each $\eta$ such that $|\eta|<\delta$, we define
\begin{align}
H^{(\eta)}(\zt) \ &:= \ H(\zt)+\eta, \notag \\
\diss_r^{(\eta)}(\zt) \ &:= \ \sum_{n=1}^NB_n(r)\left(\frac{\zt}{r}\right)^n
					+(B_0(r)+\eta), \quad\zt\in D(0;r),\notag
\end{align}
We apply Lemma \ref{L:dhull} to $\{\diss_r^{(\eta)}\}_{r\in(0;1)}$ for each $\eta$, by
taking
\[
b_n(r) \ = \ B_n(r) \quad\forall n=1,2,\dots,N, \qquad b_0(r) \ = \ B_0(r)+\eta
\]
in that lemma. For each fixed $\eta$, the family 
$\{\gr(\diss_r^{(\eta)})\}_{r\in(0,1)}$ is a continuous family of analytic discs which are 
attached to $\gr(H^{(\eta)})$ along their boundaries, and which
accumulate onto a holomorphic graph 
$\gr(\psi+\eta)$ as $r\to 1^-$. Furthermore, for $|\eta|<\delta$, each
\[
\kup^{(\eta)} \ := \ \gr(\psi+\eta)\cup 
			\left[ \ \cup_{0<r<1}\{(\zt,\diss_r^{(\eta)}(\zt))\in\CC \ | 
			\ |\zt|<r\} \ \right]\setminus\OM_1
\]
is compact. We may therefore apply the Kontinuit{\"a}tssatz to conclude that 
\begin{equation}\label{E:nbds}
U(\Delta) \ := \ \{(z,w)\in D\times\cplx \ : \ (w-\psi(z))\in\Delta\} \ = \ 
\bigcup_{\eta\in\Delta}\gr(\psi+\eta)  \ \subset \ \pi^{\OM_1}(\MOM),
\end{equation}
where $\Delta$ is any disc contained in $D(0;\delta)$
\medskip

Let ${\rm j}:\OM_1\hookrightarrow\MOM$ be the canonical holomorphic imbedding of 
$\OM_1$ into $\MOM$. As before, there exist holomorphic mappings
\[
\biss(\centerdot \ ;\eta):D(0;r)\to\MOM \ \forall r\in(0,1), \ \text{and}\quad 
\pis\centerdot \ ;\eta):D\to\MOM
\]
such that
\begin{enumerate}
\item[a)] For each $\eta$ with $|\eta|<\delta$ :
\begin{align}
\pi^{\OM_1}\circ \biss(\zt;\eta)&=(\zt,\diss_r^{(\eta)}(\zt)) \ 
\forall\zt\in D(0;r) \ \text{when} \ r\in(0,1), \ \text{and}\notag\\ 
\pi^{\OM_1}\circ \pis\zt;\eta)&=(\zt,\psi(\zt)+\eta) \ \forall\zt\in D.\notag
\end{align}
\item[b)] ${\rm j}(\zt,\diss_r^{(\eta)}(\zt))=\biss(\zt;\eta)$ wherever the left-hand
side is defined, and $\forall r\in(0,1)$.
\item[c)] ${\rm j}(\zt,\psi(\zt)+\eta)=\pis\zt;\eta)$ wherever the left-hand
side is defined.
\end{enumerate}
Arguing exactly as in the proof of Theorem \ref{T:mainC2}, we can find a disc 
$\poly\subset D(0;\delta)$, centered at the origin, and an appropriate open
set $\omega_2$ satisfying $\overline{\OM_2}\subset\omega_2\subset\OM_1$, such
that if we define 
\[
\exdt(z,w) \ := \ \begin{cases}
			\pis z;w-\psi(z)), &\text{if $(z,w)\in U(\poly)$},\\
			{} &{}\\
			{\rm j}(z,w), &\text{if $(z,w)\in\omega_2$},
			\end{cases}
\]
(where $U(\poly)$ is as defined by \eqref{E:nbds} above), then $\exdt$ is holomorphic
and well defined. Holomorphicity follows from (a) above, while (c) implies that
\begin{equation}\label{E:coincides}
\pis z;w-\psi(z)) \ = \ {\rm j}(z,\psi(z)+\{w-\psi(z)\}) \ = \ {\rm j}(z,w)\quad 
\forall(z,w)\in U(\poly)\cap\omega_2.
\end{equation}
We define $\OM^*:=U(\poly)\cup\omega_2$. As before, given any $f\in\hol(\OM_1)$, we define 
$\widetilde{f}\in\hol(\OM_*)$ by $\widetilde{f}:=\mathcal{E}(f)\circ\exdt$. In view
of \eqref{E:coincides}, $\widetilde{f}|_{\omega_2}\equiv f|_{\omega_2}$. Since
$\OM_*$ is a neighbourhood of $\gr(\psi)\cup(\bdy\times D)$, i.e. is a classical
Hartogs configuration, $\widetilde{f}$ has a holomorphic extension to $D\times D$, 
whence $f|_{\OM_2}$ has a holomorphic extension to $D\times D$. This concludes our
proof.
\end{proof}
\medskip

\begin{proof}[{\bf 3.3. The proof of Theorem \ref{T:mainCm}}] The proof of Theorem \ref{T:mainCm}
proceeds along the same lines as the proof of the first theorem. The essential difference
is that we find an $\eps>0$ such that $F(D)+\eta\subset\OM \ \forall\eta\in D(0;2\eps)^m$,
and then apply Lemma \ref{L:weird} to the pairs $(F_1,\eps),\dots,(F_m,\eps)$ to obtain
a map $G=(G_1,\dots,G_m)$ all of whose components obey the conclusions of Lemma \ref{L:weird}.
Let us write $G_j(re^{i\theta}) = \sum_{n=-N(j)}^{N(j)}B_n(r)e^{in\theta} , \ j=1,2,\dots,m$.
We apply Lemma \ref{L:1stlem} by defining
\[
b_n(r) \ := \ B_{jn}(r), \quad j\neq 0,  \qquad b_0(r) \ := \ B_{j0}(r)+\eta,
\]
for each $j=1,\dots,m$, and obtain
\begin{align}
\anal^{(\eta)}_r &:= (\anal^{(\eta)}_{1,r},\dots,\anal^{(\eta)}_{m,r}) : 
\overline{\ann}(0;r,1)\to D^m,\notag\\
\anal^{(\eta)}_r\in\hol[\ann(0;r,1)]&\cap\smoo[\overline{\ann}(0;r,1);D^m] \ \text{for
every $\eta\in D(0;\delta)^m$ and $\forall r\in(0;1)$},\notag
\end{align}
where $\delta>0$ is chosen to be so small that :
\begin{itemize}
\item $\sup_{\bdy}|H_j|+\delta<1$ for $j=1,\dots,m$;
\item $\delta+\sum_{n=-N}^N|B_{jn}(r)|/r^n < 1 \ \forall r\in(0;1]$ and $j=1,\dots,m$; and
\item $H(D)+\eta\subset\OM_1 \ \forall\eta\in D(0;\delta)^m$.
\end{itemize} 
As before, there exists a $D^m$-valued function 
$\psi:=(\psi_1,\dots,\psi_m)\in \hol(D)\cap\smoo(\overline{D};D^m)$ such
that for each $\eta\in D(0;\delta)^m$
\[
\lim_{r\to 0^+}\anal^{(\eta)}_r(\zt_0) \ = \ \psi(\zt_0)+\eta \ \text{for each
fixed $\zt_0\in D\setminus\{0\}$}.
\]
Defining
\begin{align}
H^{(\eta)} \ &:= \ (H_1(\zt)+\eta_1,\dots,H_n(\zt)+\eta_m) \quad \forall
		\eta=(\eta_1,\dots,\eta_m)\in D(0;\delta)^m, \notag\\
\kup^{(\eta)} \ &:= \ \gr(\psi+\eta)\cup
			\left[ \ \cup_{0<r<1}\gr(\anal^{(\eta)}_r) \ \right]
								\setminus\OM_1,\notag
\end{align}   
we see that properties (1)--(3) in the proof of Theorem \ref{T:mainC2} hold for
$\{\anal^{(\eta)}_r\}_{r\in(0,1)}$ and $\psi$ in our new context -- the only difference
being that the relevant functions are vector-valued, and $\eta$ varies in a
polydisc $D(0;\delta)^m$. Therefore, $\{\gr(\anal_r^{(\eta)})\}_{r\in(0,1)}$ is a 
continuous family of analytic annuli attached to $\gr(H^{(\eta)})\cup(\bdy\times D^m)$, which
accumulate onto $\gr(\psi+\eta)$ as $r\to 0^+$. The analytic annuli $\gr(\anal_r^{(\eta)})$
with $r\approx 1$ are contained in $\OM_1$. As before, the Kontinuit{\"a}tssatz tells us
that
\begin{equation}\label{E:nbdCm}
U(\Dsc) \ := \ \bigcup_{\eta\in\Dsc}\gr(\psi+\eta)  \ \subset \ \pi^{\OM_1}(\MO_1),
\end{equation}
where $\Dsc$ is any polydisc contained in $D(0;\delta)^m$.
\medskip

Arguing exactly as before, we can find a sufficiently small polydisc 
$\Dsc_*\subset D(0;\delta)^m$ centered at the origin, an appropriately chosen domain 
$\omega_2$ such that $\overline{\OM_2}\subset\omega_2\subset\OM_1$, and a mapping 
$\exdt : U(\Dsc_*)\cup\omega_2$ (here, $U(\Dsc_*)$ is as defined in \eqref{E:nbdCm} above) 
such that
\begin{itemize}
\item $\exdt\in \hol(U(\Dsc_*)\cup\omega_2)$, and 
\item $\exdt|_{\omega_2} \equiv  {\rm j}|_{\omega_2}$.
\end{itemize}
We use this mapping $\exdt$ exactly as in the previous two theorems to complete this proof.
\end{proof}
\bigskip

\section{Examples}\label{S:ex}

We begin by showing that one can construct a Hartogs-type configuration
$\gr(F)\cup(\bdy\times D)$ such that $\sup_{\overline{D}}|F|$ is as large as we
want and such that functions holomorphic in any small neighbourhood of this 
configuration extend to $D\times D$ in the manner described in Theorem \ref{T:mainC2}. 
A related {\em counterexample} to this sort of a phenomenon is the case when 
$\gr(F)$ is a Wermer disc (see Example \ref{Ex:Wermer} below for a definition). In 
that case, {\em there is no analytic continuation}, provided $\sup_{\overline{D}}|F|$
is sufficiently large. We explain under the heading Example \ref{Ex:Wermer} why this 
does not contradict Theorem \ref{T:mainC2}
\medskip

Theorem \ref{T:mainCm} is not true if $F$ occuring therein is replaced by an arbitrary
smooth, $D^m$-valued function. This is the content of Rosay's counterexample in 
\cite{rosay:crHp98}. Under Example \ref{Ex:Rosay} below, we discuss how Rosay's 
counterexample fails to meet the hypotheses of Theorem \ref{T:mainCm}.
\medskip

We begin with our first example.
\medskip

\begin{example} An example showing that given any $N\in\mathbb{N}$, there
is a $F\in\smoo^\infty(\overline{D};\cplx)$ such that $\sup_{\overline{D}}|F|>N$, 
$\sup_{\bdy}|F|<1$, and such that small connected neighbourhoods of 
$\gr(F)\cup(\bdy\times D)$ exhibit the analytic-continuation phenomenon described
in Theorem \ref{T:mainC2}.
\end{example}

We are given $N\in\mathbb{N}$. Let $\chi\in\smoo^\infty( \ [0,\infty);[0,1))$ be 
a smooth cut-off on $[0,\infty)$ such that
\begin{align}
\chi|_{[1/(N+1),\infty)} \ &\equiv \ \frac{N+(1/2)}{N+1},\notag\\ 
\chi \ &\equiv \ 0 \ \text{in a relatively open neighbourhood of $0$.}\notag
\end{align}
We define 
\[
F(re^{i\tht}) \ := \ \frac{\chi(r)}{r}e^{-i\tht}.
\]
Clearly
\begin{align}
|F(e^{i\tht}/(N+1))| \ &= \ N+(1/2) \ > \ N,\notag \\
|F(\zt)| \ &= \ \frac{N+(1/2)}{N+1} \ < \ 1 \quad\forall\zt\in\bdy,\notag \\
r|A_{-1}(r)| \ &= \ \chi(r) \ \leq \ \frac{N+(1/2)}{N+1} \ < \ 1 \quad\forall r\in(0,1].\notag
\end{align}
These conditions imply that for any connected neighbourhood $\OM_1\supset\gr(F)\cup(\bdy\times D)$,
any connected open set $\OM_2$ satisfying 
$(\bdy\times D)\subset\overline{\OM_2}\subset\OM_1\cap(\{ \ |z|\geq 1\}\times D)$, 
and for any $f\in\hol(\OM_1)$,
$f|_{\OM_2}$ extends holomorphically to $D\times D$.
\bigskip   	

Wermer presents an example of a function $g\in\smoo^\infty(\overline{D})$
\cite{hormanderWermer:uacsC^n68} with the property that $\gr(g)$ is totally 
real, but $g|_{\bdy}\equiv 0$. This allows us to define a function 
$F=Mg$ -- where $M>0$ is sufficiently large -- such that the configuration
$\gr(F)\cup(\bdy\times D)$ resists analytic continuation of the type described
in Theorem \ref{T:mainC2}. We now explain how the relevant $F$ fails to satisfy
the hypothesis of Theorem \ref{T:mainC2}.
\pagebreak

\begin{example}\label{Ex:Wermer} Wermer's disc.
\end{example}

The graph of the function
\[
g(z,\zbar) \ = \ \zbar(1-|z|^4)+i\zbar(1-|z|^2)
\]
is a totally-real surface in $\CC$. This follows from an easy computation; see 
details in \cite[Example 6.1]{hormanderWermer:uacsC^n68}. Therefore, the 
domain
\[
\dom_\delta \ := \ \{(z,w)\in\CC \ : \ |z|<1+\delta, \ |w-g(z)|<\delta \}
\]
is a pseudoconvex domain for all $\delta>0$ sufficiently small. Notice
that $\dom_\delta\supset \bdy\times D(0,\delta)$. Let $\delta^*>0$ be so small
that $D\times\{0\}\varsubsetneq \dom_{\delta^*}$ and $\dom_{\delta^*}$ is 
pseudoconvex. Then, for each domain $\OM_1\supset\gr(g)\cup(\bdy\times D(0;\delta^*))$
such that $\OM_1\subset\dom_{\delta^*}$,
there would exist a function $f\in\hol(\OM_1)$ such that $f$ does {\em not} extend
holomorphically to the bidisc $D\times D(0;\delta^*)$, because $\dom_{\delta^*}$ is a
domain of holomorphy but does not contain $D\times\{0\}$. We now define
\[
F(z) \ := \ \tfrac{1}{\delta^*}g(z), 
	\quad \DOM := \ \{(z,w) \ | \ (z,\delta^*w)\in\dom_{\delta^*}\}.
\]
By construction, $\DOM$ is a pseudoconvex domain that contains 
$\gr(F)\cup(\bdy\times D)$ but does not contain $D\times\{0\}$. By our
preceding remarks, the Hartogs-Chirka type configuration just constructed does 
not admit analytic continuation in the manner described in Theorem \ref{T:mainC2}
\medskip

Notice that $F(re^{i\tht}) = A_{-1}(r)e^{-i\tht}$, where
\[
A_{-1}(r) \ = \ \frac{r}{\delta^*}\{ \ (1-r^4) + i(1-r^2) \ \}.
\]
We will now show that
\[
r|A_{-1}(r)| \ = \ \frac{r^2}{\delta^*}\{(1-r^4)^2+(1-r^2)^2 \}^{1/2} \ \gg \ 1 
		\quad \text{\em for some $r\in(0,1] \ $,}
\]
whence Theorem \ref{T:mainC2} is inapplicable to the above configuration.
For this purpose, we will need an upper bound for the quantity $\delta^*$ introduced
above, and we make the following

\noindent{{\bf Claim :} $\delta^*<0.0061$. To see this, we refer to the 
Berndtsson-S{\l}odkowski inequality -- see \cite[Prop.2.3/(b)]{berndtsson:Lfscs93}
-- determining when a surface of the form
\[
\mathcal{S} \ = \ \{(z,w)\in\OM\times\cplx \ | \ |w-G(z)|=e^{-u(z)}\},
\]
(here $\OM$ is a domain in $\cplx$, $G$ and $u$ are smooth functions, and
$u$ is real-valued) is pseudoconvex. The desired inequality is
\begin{equation}\label{E:swallow}
\mathcal{S} \ \text{is pseudoconvex} \  \iff \ -u_{z\zbar} \ \leq \ 
	e^{2u}|G_{\zbar}|^2-e^u|G_{z\zbar}+2u_zG_{\zbar}|.
\end{equation}
For $\dom_{\delta^*}$ to be pseudoconvex, we require that the surface
$\mathcal{S}_{\delta^*}:=\{(z,w)\in D\times\cplx \ | \ |w-g(z)|=\delta^*\}$ be
pseudoconvex. Applying \eqref{E:swallow} to the surface $\mathcal{S}_{\delta^*}$
yields the following restriction on $\delta^*$.
\[
0 \ < \delta^* \ \leq \ \frac{(1-3|z|^4)^2+(1-2|z|^2)^2}{|z|\sqrt{36|z|^4+4}} 
		\quad\forall |z|\leq 1.
\]
In other words,
\[
0 \ < \ \delta^* \ \leq \ \min_{r\in[0,1]}\frac{(1-3r^4)^2+(1-2r^2)^2}{r\sqrt{36r^4+4}},
\]
and one can use any computational software package to show that the right-hand side of 
the above inequality is greater than $0.0061$. Hence the claim.
\medskip

One can also compute that 
$\max_{r\in[0,1]}r^2\{(1-r^4)^2+(1-r^2)^2\}^{1/2} \approx 0.456$. Thus
\[
\max_{r\in[0,1]}r|A_{-1}(r)| \ > \ \frac{0.456}{0.0061} \ \gg \ 1,
\]
which violates condition \eqref{E:condtnC2}.
\bigskip
 
\begin{example}\label{Ex:Rosay} Rosay's counterexample. 
\end{example}

Rosay shows that one can find an arbitrarily small, strictly pseudoconvex neighbourhood
$\OM$ of $\bdy\times D^2$ and a $D^2$-valued function $F$ such that 
$\gr(F)\cup(\bdy\times D^2)$ is holomorphically convex. Specifically 
\begin{multline}
\OM \ := \ \{(z_1,z_2,z_3)\in\cplx^3 \ : \ \{( \ |z_1|^2-1)^2+s_1|z_2|^2\}^N + 
			\left|\tfrac{z_3}{N}\right|^{2N} \\
+\alpha\{( \ |z_1|^2-1)^2+s_1|z_2|^2+|z_3|^2\} \ < \ s^{2N}+\alpha s^2\}, \notag
\end{multline}
where
\begin{itemize}
\item $s>0$ is small, and $s_1=s(1-\delta)<s$ for a fixed, small $\delta>0$, and
\item One first chooses $N$ large enough that $s_1^{2N}+1/N^{2N}<s^{2N}$, and then
chooses $\alpha$ sufficiently small so as to ensure that $\bdy\times D^2\Subset\OM$.
\end{itemize}
Write $F=(F_1,F_2)$. In Rosay's construction
\[
F_1(re^{i\theta}) \ := \ \kappa\chi(r)e^{i\theta},
\]
where $\chi\in\smoo^\infty[0,1]$ with $0\leq\chi\leq 1$, such that $\chi\equiv 1$ off a 
small relative neighbourhood of $0\in[0,1]$ and $\chi\equiv 0$ on a smaller neighbourhood
of $0$. The quantity $\kappa$ will be described presently. The function $F_2$ is 
required to be identically zero in an open set contained in $\{re^{i\theta} \ : \ r\in
\text{supp}(\chi) \ \}$, and to satisfy $\partial F_1/\partial\bar{z}_1\neq 0$ wherever
$F_1\equiv 0$. Therefore, $F_2$ will have negative Fourier modes.
\medskip

Our interest is in examining $F_1$. The constant $\kappa$ is so chosen that 
\[
\gr(F)\cap\overline{\OM} \ = \ \{(Re^{i\theta},\kappa e^{i\theta},0) \ : \ 
						\theta\in[0,2\pi) \ \},
\]
and such that $\gr(F)\cap\overline{\OM}$ is a complex-tangential curve in
the surface $\partial\OM\cap\{z_3=0\}$. Write $z_j:=x_j+iy_j, \ j=1,2$.
It is easy to determine what the magnitudes of $\kappa$ and $R$ 
(which is close to $1$) should be by visualizing 
$\omega:=\overline{\OM}\cap(\RR\times\{0\})$. Then $(R,\kappa)$ are the coordinates of 
the point of tangency, in the first quadrant, of the line through the origin 
that is tangential to $\partial\omega$ (then, the complex span of this line contains
the tangent line to the curve $\gr(F)\cap\overline\OM$). By construction, the point 
$(1,s/s_1)\in\partial\omega$ lies {\em below} the line just described, whence the line 
$x_1=x_2$ lies below this line, in the first quadrant. Thus, in the notation of Theorem
\ref{T:mainCm}
\[
\kappa/R \ :=  \ A_{11}(R)/R \ > \ 1,
\]
whence, by construction
\[
A_{11}(r)/r \ > \ 1 \quad\forall r\in\chi^{-1}\{1\}.
\]
This violates the condition \eqref{E:condtnCm}.
 

\bigskip

\begin{thebibliography}{8}

\bibitem{berndtsson:Lfscs93}
B. Berndtsson, {\em Levi-flat surfaces with circular sections}, Several Complex Variables,
Math. Notes (John Erik Fornaess ed.), Princeton University Press, Princeton, NJ (1993),
pp. 136-159.

\bibitem{bharali:sgCet00}
G. Bharali, {\em Some generalizations of Chirka's extension theorem}, Proc. Amer. Math. Soc.
{\bf 129} (2001), 3665-3669.

\bibitem{chirka:gHlnldbar?}
E.M. Chirka, {\em Generalized Hartogs' lemma and non-linear $\overline\partial$-equation},
Complex Analysis in Contemporary Mathematics (E.M. Chirka, ed.), FAZIS, Moscow (in Russian) 
(2002).

\bibitem{chirkaRosay:rpgHl98}
E.M. Chirka and J.-P. Rosay, {\em Remarks on the proof of a generalized Hartogs lemma},
Ann. Pol. Math. {\bf 70} (1998), 43-47.

\bibitem{chirkaStout:K95}
E.M. Chirka and E.L. Stout, {\em A Kontinuit{\"a}tssatz}, Topics in Complex Analysis
(P. Jac{\'o}bczak, W. Ple{\'s}niak, eds.), Banach Center Publications, Warsaw, {\bf 31} (1995),
pp. 143-150.

\bibitem{hormanderWermer:uacsC^n68}
L. H{\"o}rmander and J. Wermer, {\em Uniform approximation on compact sets in 
$\mathbb{C}^n$}, Math. Scand. {\bf 23} (1968), 5-21.

\bibitem{ivashkovichShevchishin:dncc&ems98}
S.M. Ivashkovich and V.V. Shevchishin, {\em Deformations of noncompact complex
curves, and envelopes of meromorphy of spheres} (in Russian), Mat. Sb. {\bf 189} (1998),
23-60.
  
\bibitem{rosay:crHp98}
J.-P. Rosay, {\em A counterexample related to Hartogs' phenomenon (A question by 
E. Chirka)}, Michigan Math. J. {\bf 45} (1998), 529-535.

\end{thebibliography}
\end{document}